\documentclass[11pt,oneside]{article}
\usepackage{amsfonts,amssymb,amsmath,verbatim,authblk}
\usepackage{hyperref}
\usepackage{graphicx}
\usepackage{amsbsy}
\usepackage{mathrsfs}
\newcommand{\beq}{\begin{equation}\ }
\newcommand{\eeq}{\end{equation}\ }
\begin{document}
\author[1,2]{Demetris T. Christopoulos}
\affil[1]{\small{National and Kapodistrian University of Athens, Department of Economics}}
\affil[2]{\tt{dchristop@econ.uoa.gr}, \tt{dem.christop@gmail.com} }
\title{{\itshape Polynomial regression using trapezoidal rule for computing Legendre coefficients}}
\maketitle
\begin{abstract}
We are presenting a method for computing the Fourier coefficients of a given polynomial regression by using the trapezoidal rule for numerical integration. As function basis we use the orthogonal Legendre polynomials. The results are accurate and stable compared to Forsythe's  method.\medskip
\end{abstract}
\section{Polynomial regression}
The polynomial regression technique is based on the OLS computation of the coefficients in the formal truncated series expansion of degree $m$
\beq
\label{eq:polreg}
y_{i}=\beta_0+\beta_{1}x_{i}+\beta_{2}x_{i}^2+\ldots+\beta_{m}x_{i}^m+\epsilon_{i}\,\,\epsilon_{i}\sim\,iid(0,\sigma^2)\,,\,i=1,2,\ldots,n
\eeq
or in matrix form
\beq
\label{eq-bhat}
y=\mathbf{X}\cdot\beta+\epsilon\text{ , }\epsilon|\mathbf{X}\sim\text{iid}\left(0,\sigma^2{I_{m+1}}\right)
\eeq
which has the well known OLS solution which is given by the pseudo-inverse Moore-Penrose matrix after \cite{moo-20} and \cite{pen-55}
\beq
\label{eq:ols}
\hat{\beta}=\left(\mathbf{X}^{'}\mathbf{X}\right)^{-1}\mathbf{X}^{'}y=\mathbf{X}^{+}y
\eeq
The problem is that matrix $\mathbf{X}^{'}\mathbf{X}$ can be proved, that for equidistant $x_{i}$, see \cite{for-57}, \cite{sha-97} and for increasing n is approximately the Hilbert matrix, one of the most famous ill-conditioned matrices. Early computations, see \cite{asc-58} has sown that only for a degree up to 10 we could have satisfactory outputs. Although the situation has been better now due to arbitrary precision arithmetic computations, the computational effort is still big enough. If we manage to diagonalize $\mathbf{X}^{'}\mathbf{X}$ then our task is much more easy computationally. This process has been done by Forsythe, see \cite{for-57} \& \cite{asc-58}, where a recursive method for defining orthogonal polynomials was introduced. The concept of polynomial orthogonality there had the sense of discrete orthogonality, i.e. two polynomial $\phi_{k}(x),\phi_{l}(x)$ are said to be orthogonal over a set of abscissae $\left\{x_1,x_2,\ldots,x_n\right\}$ if the next vanishing equation holds
\beq
\label{eq-polorth1}
\sum_{i=1}^{n}{\phi_{k}(x_i)\phi_{l}(x_i) }=0
\eeq
This is nothing else than the zero common Euclidean real inner product
\beq
\label{eq-polorth2}
\left\langle \boldsymbol{\phi_{k}},\boldsymbol{\phi_{l}}\right\rangle= \boldsymbol{\phi_{k}}^{'}\,\boldsymbol{\phi_{l}}=0
\eeq
By using this procedure instead of directly computing $\beta$ coefficients of \ref{eq:polreg} we compute the coefficients of next truncated series
\beq
\label{eq:orthopol1}
y_{i}=\beta_0\phi_{0}(x_{i})+\beta_{1}\phi_{1}(x_{i})+\beta_{2}\phi_{2}(x_{i})+\ldots+\beta_{m}\phi_{m}(x_{i})+\epsilon_{i},\,\epsilon_{i}\sim\,iid(0,\sigma^2)
\eeq
or in matrix form
\begin{equation}
\label{eq:orthopol2}
y=\mathbf{\Phi}\cdot\beta+\epsilon\text{ , }\epsilon|\mathbf{\Phi}\sim\text{iid}\left(0,\sigma^2{I_{m+1}}\right)
\end{equation}

where the design matrix is
\begin{equation}
\label{eq:orthopol3}
\mathbf{\Phi}
=
\begin{bmatrix}
\boldsymbol{\phi_{0}}&\boldsymbol{\phi_{1}}&\ldots&\boldsymbol{\phi_{m}}
\end{bmatrix}
=
\begin{bmatrix}
1&\phi_{1}(x_1)&\ldots&\phi_m(x_1)\\
1&\phi_{1}(x_2)&\ldots&\phi_m(x_2)\\
\vdots&\vdots&\ldots&\vdots\\
1&\phi_{1}(x_n)&\ldots&\phi_m(x_n)\\
\end{bmatrix}
\end{equation}
Our coefficients are now simply the well known Fourier coefficients
\beq
\label{eq:four1}
\beta_{j}=\frac{\left\langle y,\boldsymbol{\phi_{j}}\right\rangle}{\left\langle \boldsymbol{\phi_{j}},\boldsymbol{\phi_{j}}\right\rangle}=
\frac{\sum_{i=1}^{n}{y_i\,\phi_{j}(x_i) }}{\sum_{i=1}^{n}{\phi_{j}(x_i)^2 }}
\eeq
If we take orthonormal polynomials, see \cite{hav-92}, i.e. if it holds that
\beq
\label{eq-orthnorm1}
\left\langle \boldsymbol{\phi_{k}},\boldsymbol{\phi_{l}}\right\rangle= \boldsymbol{\phi_{k}}^{'}\,\boldsymbol{\phi_{l}}=\delta_{kl}= \begin{cases}1 & k=l\\ 0 & k\neq l\end{cases}
\eeq
then $\mathbf{\Phi}^{'}\mathbf{\Phi}=I_{m+1}$ and the coefficients are simply
\beq
\label{eq:four2}
\beta_{j}=\mathbf{\Phi}^{'}y=\left\langle y,\boldsymbol{\phi_{j}}\right\rangle=
\sum_{i=1}^{n}{y_i\,\phi_{j}(x_i) }\,,\,j=0,\ldots,m
\eeq
The concept of orthogonality or orthonormality is a linear algebraic term and is independent of the chosen function basis representation. It is an elementary exercise, see \cite{spi-74} page 61, that starting from the linearly independent set of monomials $\left\{x^i,i=0,1,2,\ldots,m,\forall m\in\aleph\right\}$ and by using Gram-Schmidt orthonormalisation process we can end to the normalised Legendre polynomials
\begin{equation}
\label{eq:normleg}
\mathscr {P}_k(x)=\sqrt{m+\frac{1}{2}}\,\,P_{k}(x)=
\sqrt{m+\frac{1}{2}}\,\,\frac{1}{2^k\,k!}\,{\frac {d^{k}}{d{x}^{k}}}(x^2-1)^k  
\end{equation}

The above polynomials are orthonormal in the interval $[-1,1]$
\begin{equation}
\label{eq:leg11}
\int_{-1}^{1}{\mathscr {P}_k(x)\mathscr {P}_l(x)}=\delta_{kl}
\end{equation}
Now we can expand every function in a Legendre series expansion
\beq
\label{eq:legser}
f(x)=\sum_{i=1}^{m}{c_k\,\mathscr {P}_k(x)dx}
\eeq
with the Fourier coefficients given by
\begin{equation}
\label{eq:coeff1}
\beta_k=\left\langle f(x) , \mathscr {P}_{k}(x) \right\rangle
=\int_{-1}^{1}{f(x)\mathscr {P}_{k}(x)\,dx}
\end{equation}
If we follow the guides of \cite{sha-97} and make the linear transformation 
\beq
T_{2}:[a,b]\rightarrow[-1,1]\quad,\quad{T_{2}(x)=\frac{2x-a-b}{b-a}}
\label{eq:t21}
\eeq
in order to convert our initial range $[a,b]$ to the $[-1,1]$, where many orthogonal polynomials are defined, then the Forsythe polynomials are just a scale version of Legendre polynomials. So, the norm we have chosen does not play any other role except for the simplicity of computations. If we choose the $l_2$ norm, then we can proceed like Forsythe and construct a set of orthogonal polynomials for solving our polynomial regression problem. \\
The discrete case is
\begin{equation}
\label{eq:orthonorm1}
y_{i}=f(x_i)=\beta_{0}\mathscr {P}_{0}(x_{i})+\beta_{1}\mathscr {P}_{1}(x_{i})+\beta_{2}\mathscr {P}_{2}(x_{i})+\ldots+\beta_{m}\mathscr {P}_{m}(x_{i})\\
\end{equation}
It is obvious to think about computing the Fourier coefficients \ref{eq:coeff1} by a numerical approximation of the relevant integral.
For the equidistant case with $x_{j+1}-x_{j}=1$ we have that
\begin{equation}
\label{eq:coeff2}
\beta_k=\left\langle f(x) , \mathscr {P}_{k}(x) \right\rangle
\approx{\sum_{i=1}^{n}{y_i\,\mathscr {P}_{k}(x_i)}}
\end{equation}	
Now we have approximated the integral via the orthogonal rule. We can also use trapezoidal method in order to increase the accuracy. By comparing \ref{eq:coeff2} and \ref{eq:four2} we see that the latter is just the left orthogonal Riemannian approximation for the continuous case \ref{eq:coeff1}. If we had use the simple Legendre polynomials $P_k(x)$ then our coefficients could be
\begin{equation}
\label{eq:coeff3}
\beta_k=
\frac{\left\langle f(x) , {P}_{k}(x) \right\rangle}
{\left\langle {P}_{k}(x) , {P}_{k}(x) \right\rangle}
\end{equation}
For the equidistant case $x_{j+1}-x{j}=h$ we have that
\begin{equation}
\label{eq:coeff4}
\beta_k
\approx{\frac{\sum_{i=1}^{n}{y_i\,{P}_{k}(x_i)\,h}}
{\sum_{i=1}^{n}{{P_{k}}^2(x_i)\,h}}}=
\frac{\sum_{i=1}^{n}{y_i\,{P}_{k}(x_i)}}
{\sum_{i=1}^{n}{{P_{k}}^2(x_i)}}
\end{equation}
which is just \ref{eq:four1} for Fourier coefficients.

\section{Trapezoidal estimation of Fourier coefficients}
Our task is to compute the integrals of our Fourier coefficients, \ref{eq:coeff2} for normalised or \ref{eq:coeff3} for simple orthogonal polynomials by using the trapezoidal rule of numerical integration. We shall constraint in the equidistant case, since we have closed formulas using less arithmetic operations. \\
For the case of simple orthogonal polynomials we have that 
\begin{equation}
\label{eq:coeff5}
\beta_k
\approx{\frac{y_{1}\,P_k(x_1)+2\,\sum_{i=2}^{n-1}{y_i\,{P}_{k}(x_i)}+y_{n}\,P_k(x_n)}
{{P_{k}}^2(x_1)+2\,\sum_{i=1}^{n}{{P_{k}}^2(x_i)\,h}+{P_{k}}^2(x_n)}}
\end{equation}
If we use orthonormal polynomials we have the estimation
\begin{equation}
\label{eq:coeff6}
\beta_k
\approx{\frac{h}{2}\,\left(y_1\,\mathscr {P}_{k}(x_1)+2\,\sum_{i=2}^{n-1}{y_i\,\mathscr {P}_{k}(x_i)}+y_n\,\mathscr {P}_{k}(x_n)\right) }
\end{equation}	
The total sum of squares is almost identical for the two cases and for the simple OLS regression by mean of \ref{eq:ols} with $\mathbf{X}=\mathbf{\Phi}$.
\section{A numerical example}
\begin{figure}
\begin{center}
\caption{The plot of the function used at numerical example}\label{fig:01}
\vspace{0.5in}  \includegraphics[width=6cm,height=6cm]{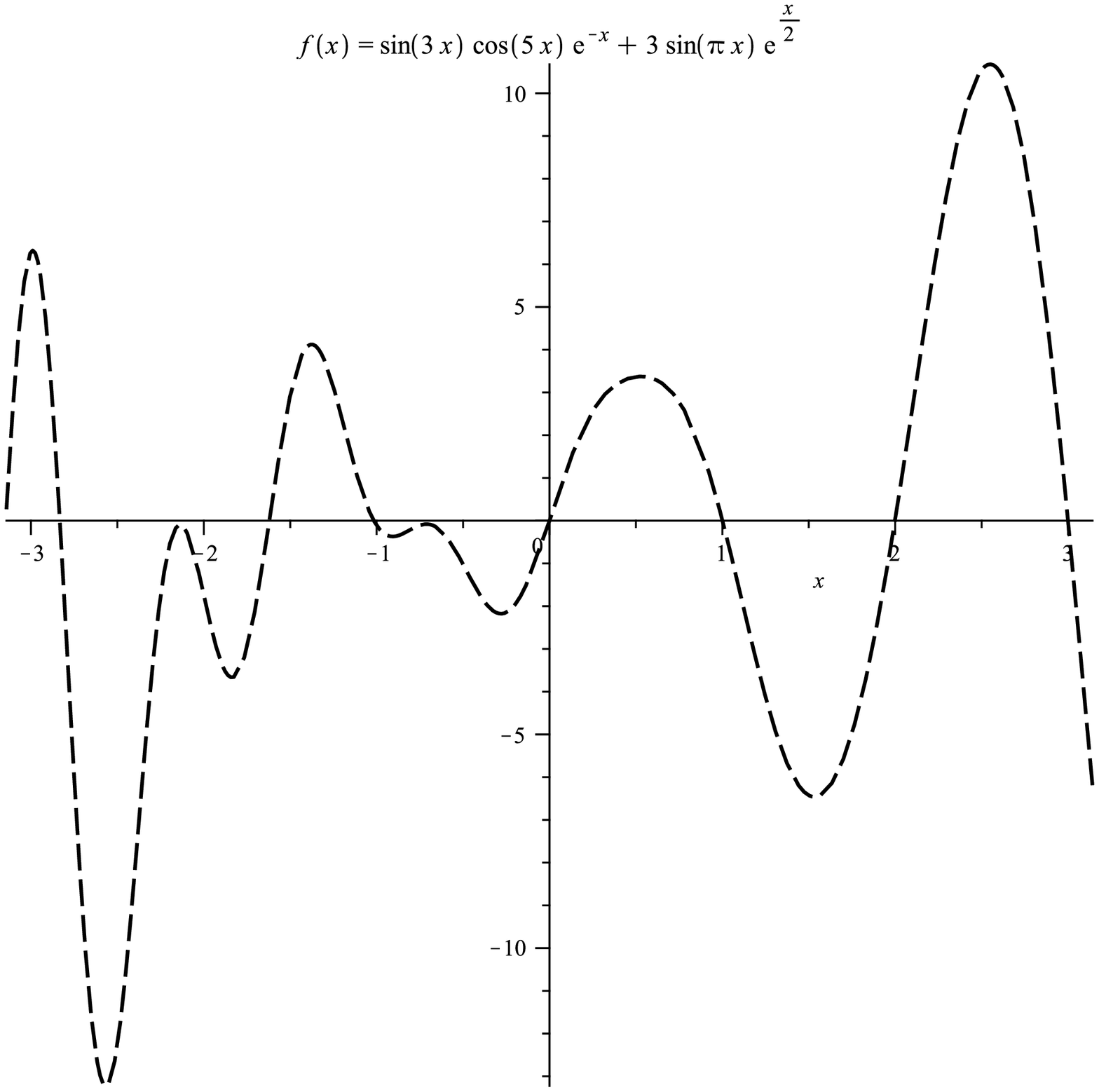}
\vspace{0.5in}
\end{center}
\end{figure}

\begin{figure}
\begin{center}
\caption{Legendre series approximation}\label{fig:02}
\vspace{0.5in}  \includegraphics[width=6cm,height=6cm]{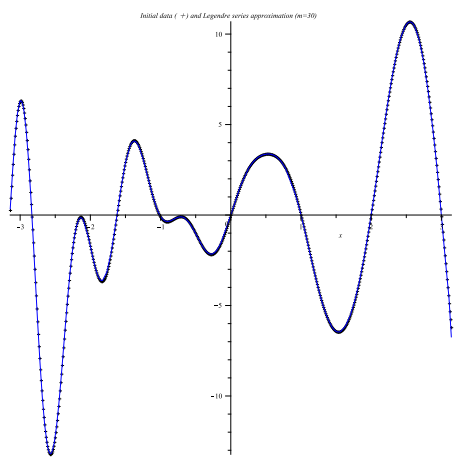}
\vspace{0.5in}
\end{center}
\end{figure}

Let us consider the known function:
\beq
{f:}\left[-\pi,\pi\right]\rightarrow\Re,\,\,f \left( x \right) =\sin \left( 3\,x \right) \cos \left( 5\,x \right)
{e^{-x}}+3\,\sin \left( \pi \,x \right) {e^{\frac{x}{2}}} 
\label{eq:ftest}
\eeq
at an equal spaced grid $x_{i},i=0,\ldots,628$. The graph of the function is presented in Figure \ref{fig:01}. This function is a smooth function, $f\in{C^{\infty}}$ , it has 6 local maxima and 5 local minima inside the interval $[-\pi,\pi]$. It is a rather complicated function, for example it has one local minimum and maximum in the small interval $[-1,-0.5]$, so the task of recovering this shape is difficult.\\
By using floating point arithmetic with 32 digits of accuracy we can obtain the next Taylor polynomial of $30^{th}$ degree\\
$$
\begin{array}{lll}
T_{{30}} \left( x \right)&=&+   0.0000016359\,{x}^{30}- 0.0000097947\,{x}^{29}- 0.000013156\,{x}^{28}\\
&&+   0.0001336926\,{x}^{27}+ 0.0000419458\,{x}^{26}- 0.001477436\,{x}^{25}\\
&&+   0.0007170316\,{x}^{24}+ 0.0131083743\,{x}^{23}- 0.015365948\,{x}^{22}\\
&&-   0.091642089\,{x}^{21}+  0.168431320\,{x}^{20}+  0.488498913\,{x}^{19}\\
&&-   1.270259183\,{x}^{18}-  1.866729762\,{x}^{17}+  6.956422865\,{x}^{16}\\
&&+   4.386886795\,{x}^{15}- 27.70967780\,{x}^{14}-   2.236993754\,{x}^{13}\\
&&+  78.47606337\,{x}^{12} - 23.60696807\,{x}^{11}- 151.3009259\,{x}^{10}\\
&&+  86.58066121\,{x}^{9}+  184.9106325\,{x}^{8}-   142.4655789\,{x}^{7}\\
&&- 125.9202816\,{x}^{6}+   121.2621435\,{x}^{5}+    33.94478037\,{x}^{4}\\
&&-  54.82504110\,{x}^{3}+    1.71238898\,{x}^{2}+   12.42477796\,x\\
\end{array}
$$
Although the above polynomial is not identical to the initial function outside approximately the interval $[-1.5,1.5]$ it is a representation that carries a lot of information about the function since it can give the derivatives until the $30^{th}$ order. So, our task is (i) to recover as many as possible coefficients of the above series expansion and (ii) to approximate the functional data with the smallest possible error.\\
In order to avoid multicollinearity problems due to lower accuracy we are using 32 digits in our arithmetic operations and we are transforming to the interval $[-1,1]$ both $x_i\,\&\,y_i$ data. After finishing our coefficient computations we are performing the inverse $T_2$-transform and return to our initial data scale.\\
Results are presented at Table \ref{tab:1} while the sum of squares for both cases, the transformed to $[-1,1]$ and the initial, are given at Table \ref{tab:2}.\\
The relevant plots of all Legendre series are indistinguishable from the original data, see Figure \ref{fig:02}.\\

As a benchmark to our effort we shall compare our results with those obtained by using \cite{for-57} method as has been implemented in FORTRAN 90 by \cite{jpm-13}. We find that under double precision arithmetic, i.e. with 16 digits accuracy, the solution divergences very fast from the true series expansion. The \cite{for-57} polynomial coefficients  after inverse transforming to the initial domain are presented at Table \ref{tab:3}. 

\section{Discussion}
The times for computing the coefficients were
 $\left(P_k(x),\mathscr {P}_{k}(x),{P_k}^{OLS}(x)\right)=(16.895,11.013,24.820)$ CPU seconds in a typical Intel Core i5 CPU with 4 GB RAM memory and by using Maple program. We observe that the use of orthonormal polynomials is reducing the computational time.\\
If we decrease our accuracy to 16 digits in order to be compatible with FORTRAN we obtain similar results, see Table \ref{tab:4}. Thus our methods still found converged and suitable outputs compared to the \cite{for-57} orthogonal polynomial method.

\newpage

\begin{table}
\caption{Legendre coefficients}\label{tab:1}
$$
\begin {matrix}
k&P_k(x)&\mathscr {P}_{k}(x)&{P_k}^{OLS}(x)&T_{k}(x) \\
 0&- 0.0001255282& 0.001814486&- 0.002035062& 0.0\\ 1& 12.33649& 12.44636& 12.50049&
 12.42478\\ 2& 1.852769& 1.629086& 1.821539&
 1.712389\\ 3&- 54.36178&- 55.30185&- 56.14278&-
 54.82504\\ 4& 32.78326& 34.56230& 32.97494&
 33.94478\\ 5& 121.7605& 124.2576& 128.1177&
 121.2621\\ 6&- 122.8261&- 127.6625&- 122.5218&-
 125.9203\\ 7&- 148.0910&- 151.0391&- 159.2726&-
 142.4656\\ 8& 180.7955& 187.3529& 178.6751&
 184.9106\\ 9& 98.54291& 100.3063& 110.2045&
 86.58066\\ 10&- 147.9384&- 153.1944&- 144.4110&-
 151.3009\\ 11&- 36.87822&- 37.37939&- 44.82297&-
 23.60697\\ 12& 76.54964& 79.28914& 73.50992&
 78.47606\\ 13& 7.061322& 7.061628& 10.78721&-
 2.236994\\ 14&- 26.87537&- 27.85941&- 25.26188&-
 27.70968\\ 15&- 0.09374941&- 0.04223660&- 1.330764&
 4.386887\\ 16& 6.680829& 6.933201& 6.111261&
 6.956423\\ 17&- 0.3086695&- 0.3285492&- 0.01411589&
- 1.866730\\ 18&- 1.204759&- 1.251863&- 1.065779&-
 1.270259\\ 19& 0.08555660& 0.08970893& 0.03517729&
 0.4884989\\ 20& 0.1590441& 0.1654770& 0.1352121&
 0.1684313\\ 21&- 0.01253891&- 0.01309178&-
 0.006412513&- 0.09164209\\ 22&- 0.01527912&-
 0.01591658&- 0.01241044&- 0.01536595\\ 23&
 0.001141524& 0.001189852& 0.0006253121& 0.01310837
\\ 24& 0.001042632& 0.001087308& 0.0008050851&
 0.0007170316\\ 25&- 0.00006491531&- 0.00006761426&-
 0.00003629592&- 0.001477436\\ 26&- 0.00004797304&-
 0.00005007409&- 0.00003507896& 0.00004194580\\ 27&
 0.000002124622& 0.000002212277& 0.000001186322& 0.0001336926
\\ 28& 0.000001335500& 0.000001395001&
 0.0000009223072&- 0.00001315592\\ 29&-
 0.00000003068878&- 0.00000003195140&- 0.00000001691387&-
 0.000009794687\\ 30&- 0.00000001699460&-
 0.00000001776155&- 0.00000001106804& 0.000001635893
 \end {matrix}
$$
\end{table}

\begin{table}
\caption{Total Sum of Squares}\label{tab:2}
$$
\begin {matrix}
Domain&P_k(x)&\mathscr {P}_{k}(x)&{P_k}^{OLS}(x)\\
[-1,1]\times [-1,1]&218.869& 219.183& 219.095\\
[-\pi,\pi]\times[y_{min},y_{max}]&13132.6& 13177.4& 13169.9
\end{matrix}
$$
\end{table}

\begin{table}
\caption{Forsythe orthogonal polynomial coefficients}\label{tab:3}
$$
\begin {matrix}
k&{P_k}^{Forsythe}(x)&T_{k}(x) \\
0& 0.05094099& 0.0\\ 1& 0.9398571& 12.42478\\ 2& 0.1309315& 1.712389
\\ 3& 0.1168369&- 54.82504\\ 4&
 0.1049183& 33.94478\\ 5&- 0.2143114& 121.2621
\\ 6&- 0.04026457&- 125.9203\\ 7&-
 0.08311646&- 142.4656\\ 8&- 0.006660081& 184.9106
\\ 9& 0.02541863& 86.58066\\ 10&
 0.001024787&- 151.3009\\ 11& 0.002504286&- 23.60697
\\ 12&- 0.003783044& 78.47606\\ 13
& 0.001420267&- 2.236994\\ 14&- 0.0002494057&-
 27.70968\\ 15&- 0.0002088982& 4.386887
\\ 16& 0.0003358103& 6.956423\\ 17
&- 0.0002421829&- 1.866730\\ 18& 0.00008259028&-
 1.270259\\ 19& 0.00001781517& 0.4884989
\\ 20&- 0.00004319505& 0.1684313
\\ 21& 0.00002519070&- 0.09164209
\\ 22& 0.00000009195903&- 0.01536595
\\ 23&- 0.000007925483& 0.01310837
\\ 24& 0.000002842523& 0.0007170316
\\ 25& 0.0000009651675&- 0.001477436
\\ 26&- 0.0000007763782& 0.00004194580
\\ 27&- 0.000000005683953& 0.0001336926
\\ 28& 0.0000001149406&- 0.00001315592
\\ 29&- 0.00000001691387&- 0.000009794687
\\ 30&- 0.00000001106804& 0.000001635893
 \end {matrix}
$$
\end{table}

\begin{table}
\caption{Legendre coefficients for 16 digits accuracy}\label{tab:4}
$$
\begin {matrix}
k&P_k(x)&\mathscr {P}_{k}(x)&{P_k}^{OLS}(x)&T_{k}(x) \\
0&- 0.0001255282& 0.001814486&- 0.002035063& 0.0\\ 1& 12.33649& 12.44636& 12.50049&
 12.42478\\ 2& 1.852769& 1.629086& 1.821539&
 1.712389\\ 3&- 54.36178&- 55.30185&- 56.14278&-
 54.82504\\ 4& 32.78326& 34.56230& 32.97494&
 33.94478\\ 5& 121.7605& 124.2576& 128.1177&
 121.2621\\ 6&- 122.8261&- 127.6625&- 122.5218&-
 125.9203\\ 7&- 148.0910&- 151.0391&- 159.2726&-
 142.4656\\ 8& 180.7955& 187.3529& 178.6751&
 184.9106\\ 9& 98.54291& 100.3063& 110.2045&
 86.58066\\ 10&- 147.9384&- 153.1944&- 144.4110&-
 151.3009\\ 11&- 36.87822&- 37.37939&- 44.82297&-
 23.60697\\ 12& 76.54964& 79.28914& 73.50992&
 78.47606\\ 13& 7.061322& 7.061628& 10.78721&-
 2.236994\\ 14&- 26.87537&- 27.85941&- 25.26188&-
 27.70968\\ 15&- 0.09374941&- 0.04223660&- 1.330764&
 4.386887\\ 16& 6.680829& 6.933201& 6.111261&
 6.956423\\ 17&- 0.3086695&- 0.3285492&- 0.01411583&
- 1.866730\\ 18&- 1.204759&- 1.251863&- 1.065779&-
 1.270259\\ 19& 0.08555660& 0.08970893& 0.03517728&
 0.4884989\\ 20& 0.1590441& 0.1654770& 0.1352121&
 0.1684313\\ 21&- 0.01253891&- 0.01309178&-
 0.006412512&- 0.09164209\\ 22&- 0.01527912&-
 0.01591658&- 0.01241044&- 0.01536595\\ 23&
 0.001141524& 0.001189852& 0.0006253120& 0.01310837
\\ 24& 0.001042632& 0.001087308& 0.0008050850&
 0.0007170316\\ 25&- 0.00006491531&- 0.00006761426&-
 0.00003629592&- 0.001477436\\ 26&- 0.00004797304&-
 0.00005007409&- 0.00003507895& 0.00004194580\\ 27&
 0.000002124622& 0.000002212277& 0.000001186322& 0.0001336926
\\ 28& 0.000001335500& 0.000001395001&
 0.0000009223070&- 0.00001315592\\ 29&-
 0.00000003068878&- 0.00000003195140&- 0.00000001691387&-
 0.000009794687\\ 
 30&- 0.00000001699460&- 0.00000001776155&- 0.00000001106803& 0.000001635893
 \end{matrix}
$$
\end{table}

\end{document}